\documentclass[a4paper,11pt]{article} 
\usepackage{amsmath,amssymb,latexsym}
\catcode`\é = \active\def é{\'e} \catcode`\è = \active\def è{\`e}
\catcode`\ç = \active\def ç{\c{c}} \catcode`\à = \active\def
à{\`a} \catcode`\ê = \active \def ê{\^e} \catcode`\ù = \active
\def ù{\`u} \catcode`\ô = \active \def ô{\^o} \catcode`\û =
\active \def û{\^u} \catcode`\î = \active \def î{\^{\i}}
\catcode`\â = \active \def â{\^a} \catcode`\ö = \active \def
ö{{\"o}}

\newtheorem{thm}{Theorem}[section]  \newtheorem{cor}[thm]{Corollary}
\newtheorem{lem}[thm]{Lemma}   
  \newtheorem{prop}[thm]{Proposition}

\def\remark{\refstepcounter{thm}\bigskip\noindent\bf Remark \thethm\rm\ }
\newcommand{\preuve}[1][\!\!]{\bigskip\noindent{\bf Proof #1. \ \ }}
\def\fin{\hfill$\Box$\\}
\def\ititem#1{\bigskip \par \noindent \it #1 \rm}
\def\aaa{{\cal A}}\def\cc{{\cal C}}\def\dd{{\cal D}}
\def\lll{{\cal L}}
\def\kk{{\cal K}}
\def\mm{{\cal M}}
\def\ss{{\cal S}}

\def\R{\mathbb R}
\def\H{\mathbb H}     
\def\D{\partial}\def\eps{\varepsilon}\def\phi{\varphi}
\def\norm#1{\left\Vert#1\right\Vert}

\def\set#1{\left\{#1\right\}}
\def\seq#1{\left<#1\right>}
\def\sep#1{\left(#1\right)}\def\adf#1{\left[#1\right]}
\def\Re{{\mathrm Re\,}} 
\def\defegal{\stackrel{\text{\rm def}}{=}}

\def\Vnl{V_{\text{\rm\scriptsize{nl}}}}
\def\Ve{V_{\text{\rm\scriptsize{e}}}}
\def\V{V}

\def\Vinf{V_{\infty}}

\def\mminf{{\cal M}_\infty}

\def\Bperp{{\cal B}_{\perp}}

\def\Binf{{\cal B}_\infty}
\def\Binfp{{\cal B}_{\infty,\perp}}

\def\Kinf{K_\infty}

\def\ai{\alpha_\infty}
\def\Ai{A_\infty}

\def\Vd{V_{\text{\rm\scriptsize{diff}}}}
\def\entoure#1{\framebox{\scriptsize#1}}

\def\Hess{\text{\rm{Hess}}}
\def\Id{\textrm{Id}}

\def\Hun{{\rm(H1) }}
\def\Hdeux{{\rm (H2) }}
\def\Hdeuxbis{{\rm (H2bis) }}
\def\Htrois{{\rm (H3) }}

\def\Section#1{\section{#1}}\def\Subsection#1{\subsection{#1}}

\def\coord{Laboratoire de Mathématiques, UFR Sciences exactes et
naturelles, Université de Reims-Champagne Ardennes, Moulin de la
Housse, BP 1039 51687 Reims cedex 9, \tt herau@univ-reims.fr \rm }
\date{January 30, 2006}
\title{Short and long time behavior of the Fokker-Planck equation
in a confining potential and applications}
\author{Fr{\'e}d{\'e}ric H{\'e}rau \footnote{\coord} \\
 Universit{\'e} de Reims }

\begin{document}
 \maketitle \thispagestyle{empty}
 \bibliographystyle{plain}
 \begin{quote}
   \bf Abstract\rm: We consider the linear Fokker-Planck equation in
   a  confining potential in space dimension $d \geq 3$.
   Using spectral methods, we
   prove bounds on the derivatives of the
   solution for short and long time, and give some applications.
 \end{quote}

\begin{quote}
   \bf Résumé\rm : On considère l'équation de Fokker-Planck avec
   un potentiel confinant
   en dimension $d\geq 3$. Avec des méthodes spectrales
   on donne des bornes sur les dérivées de la solution en temps
   petit et grand, et quelques applications.
 \end{quote}


\Section{Introduction and results}

In this article, we consider the linear Fokker-Planck equation in
$\R_{x,v}^{2d}$ for $d\geq 3$ which reads after scaling

\begin{equation}
  \label{VPFPS}
  \left\{
  \begin{aligned}
   &  \D_{t}f+v.\D_{x}f - \D_x \V.\D_{v}f
      -\gamma \D_{v}.\left(\D_{v}+v\right)f=0, \\
    &   f|_{t=0} = f_0,
    \end{aligned}
    \right.
\end{equation}
where $\V $ is a given external confining potential, $\gamma$ is a
positive physical constant, and $f$ is the distribution function
of the particules. This equation is a linear model for plasmas or
stellar systems, and  $\gamma$ has to be understood as a
friction-diffusion coefficient. The aim of this article is to
study the short and long time behavior of the solution of this
equation, without the help of the explicit Green function, which
is known only in very special cases (i.e. $\V $ quadratic), and
give an application to a mollified Vlasov-Poisson-Fokker-Planck
equation.

\bigskip Let us now precise our notations and hypothesis. For the
potential $V$, we suppose  the following:
$$ \Hun \ \ \ \ \  \ \ e^{-\V} \in
\ss(\R^d_x), \ \ \ \text{with}  \ \ \ \V \geq 0 \ \  \textrm{ and
} \ \   \V'' \in W^{\infty,\infty}.
$$
 Note that the assumption $0 \leq \V $ can be
relaxed by assuming that $V$ is bounded from below and adding to
it a sufficiently large constant. Let us also note that these
assumptions easily imply that $\V \in \cc^\infty(\R^d_x)$ and  $
\lim_{x \longrightarrow \infty} \V(x) = + \infty$. We introduce
now the so-called Maxwellian, which is the unique $L^1$-normalized
steady solution of equation (\ref{VPFPS}):
\begin{equation} \label{defmme}
\mm(x,v) = \frac{e^{-(v^2/2 + \V(x))}}{\int e^{-(v^2/2 +
\V(x))}dxdv}.
\end{equation}
 To this function we associate  a
 weighted space
 here called  $B^2$ built from the standard $L^2$
space after conjugation with a half power of the Maxwellian:
 \begin{equation} \label{defb2}
B^2 \defegal \mm^{1/2} L^2 =  \set{ f \in \dd' \text{ s. t. }
f/\mm \in L^2(\mm dxdv)}
\end{equation}
with the natural norm defined by
\begin{equation*}
\begin{split}
& \norm{f}_{B^2}^2=\int \sep{{f}/{\mm}}^2 \mm dxdv.
\end{split}
\end{equation*}
This space is standard (e.g. \cite{HN04},\cite{HelN04}) for the
study of the Fokker Planck operator
\begin{equation} \label{Ke}
K = v.\D_{x} -\D_x \V.\D_{v}
      -\gamma \D_{v}.\left(\D_{v}+v\right).
\end{equation}
for which $\mm$ is the unique fundamental state. It is shown in
\cite{HelN04} that $K$ is maximal accretive with $\mm^{1/2} \ss$
as a core.  It is also shown in \cite{HN04} that the associated
semi-group has smoothing properties under slightly weaker
hypotheses. Anyway a number of results from there are still true
(see Section \ref{functional} here and in particular Remark
\ref{difference}). The first Theorem concerns accurate estimates
about the short time behavior of the associated semi-group:

\begin{thm} \label{disp0}
There exists a constant $C$  such that
 for all $ t >0 $, we have the following: \\
 i)  $(-\D_v + v)e^{-tK}  \textrm{ is bounded by } C (1+t^{-1/2})
$ and\\
ii) $ (-\D_x + \D_x V) e^{-tK}  \textrm{ is bounded by } C
(1+t^{-3/2}),
$ \\
as bounded operators on $B^2$. Here  $C$ depends only
 on $\norm{V''}_{L^\infty}$ (and $\gamma$).
\end{thm}

In order to study the long time behavior  of the system
(\ref{VPFPS}) we  introduce an additional hypothesis on $\V$. We
first define an intermediate operator called the Witten Laplacian
(on $0$-forms) naturally associated  to the linear Fokker-Planck
operator $K$
\begin{equation} \label{lambrine}
\Lambda^2 =-\gamma \D_v(\D_v +  v) -\gamma \D_x(\D_x + \D_x\V ).
\end{equation}
The closure of this operator defined in $B^2$ has also $0$ as
single eigenvalue for the eigenfunction $\mm$. We shall assume
the following:
\begin{equation*}
 \Hdeux \ \ \ \ \ \ \ \  \left\{
 \begin{aligned}
& \textrm{Operator $\Lambda^2$ has a spectral gap
in $B^2$} \\
&  \textrm{with first non-zero eigenvalue denoted $\alpha$}.
\end{aligned}
\right.
\end{equation*}

This hypothesis may seem complicated, but in the particular case
when $\D_x \V \longrightarrow \infty$  it is immediate since then
$\Lambda^2$ is with compact resolvent. Under this hypothesis, we
have the following result about the (short and) long time behavior
of the solution of the Fokker-Planck equation:

\begin{thm} \label{mainn}
Suppose conditions \Hun and \Hdeux are fulfilled. Then there
exists constants $C$ and $A$  depending only
 on $\norm{V''}_{L^\infty}$ (and $\gamma$) such that if
   $f(t)$ is the solution of
(\ref{VPFPS}) for an $L^1$-normalized initial datum
$f_0 \in B^2$, we have \\
i) $\norm{f(t) -\mm}_{B^2} \leq 3 e^{-\alpha t/A}
\norm{f_0 -\mm}_{B^2}$, \\
ii) $\norm{(-\D_v +v)f(t)}_{B^2} \leq C(1+t^{-1/2})
 e^{-\alpha t/A} \norm{f_0 }_{B^2}$, \\
iii) $\norm{(-\D_x +\D_xV)f(t)}_{B^2} \leq C(1+t^{-3/2})
e^{-\alpha t/A} \norm{f_0 }_{B^2}$.
\end{thm}

We give now an application to a non-linear problem. We want to
study the following mollified Vlasov-Poisson-Fokker-Planck
equation
\begin{equation}
  \label{VPFPAintro}
  \left\{
  \begin{aligned}
   &  \D_{t}f+v.\D_{x}f -(E+ \D_x \V).\D_{v}f
      -\gamma \D_{v}.\left(\D_{v}+v\right)f=0 \\
   & E(t,x) \defegal \D_x \Vnl(t,x)
   = -\zeta * \frac{\kappa}{|S^{d-1}|}
    \frac{x}{|x|^d}*_x \rho(t,x)  \ \
     \ \ \
     \text{where } \rho(t,x) = \int     f(t,x,v)dv \\
    &   f|_{t=0} = f_0
    \end{aligned}
    \right.
\end{equation}
where $\zeta \in \ss$ (depending only on $x$). Here $\kappa \in
\R$ has to be understood as the total charge of the system. In the
usual VPFP equation there is no convolution with $\zeta$, but we
were not able to reach similar result in this case. The unique
steady state of this equation is given by
$$
\mminf(x,v) = \frac{e^{-(v^2/2 + \V(x)+\Vinf(x))}}{\int e^{-(v^2/2
+ \V(x)+\Vinf(x))}dxdv},
$$
where $\Vinf$ is a solution of the following Poisson-Emden type
equation
\begin{equation} \label{VVinf}
-\Delta \Vinf = \kappa \zeta *_x \frac{ e^{-(V+\Vinf)}}{\int
e^{-(V(x)+\Vinf(x))}dx}.
\end{equation}
It is easy to see that under hypothesis \Hun this equation has a
unique (Green) solution $\Vinf \in W^{\infty,\infty}$ thanks to
the ellipticity properties of the Laplacian. We immediately check
that the associated total potential $\V+\Vinf$ satisfies
hypothesis (H1), and that $\mminf \in \ss$.
 We define also
the associated spaces $\Binf^2 =  \mminf^{1/2} L^2$ and we impose
in addition that $V+\Vinf$ satisfies an hypothesis of type (H2):
 As in (\ref{lambrine}) we define the corresponding  Witten Laplacian
\begin{equation*}
\Lambda_\infty^2 =-\gamma \D_v(\D_v +  v) -\gamma \D_x(\D_x +
\D_xV +\D_x \Vinf )
\end{equation*}
which closure  in $\Binf^2$ has $0$ as single eigenvalue
associated with the eigenfunction $\mminf$. We shall assume  the
following:
\begin{equation*}
 \Hdeuxbis \ \ \ \ \ \ \ \  \left\{
 \begin{aligned}
& \textrm{Operator $\Lambda_\infty^2$ has a spectral gap
in $\Binf^2$} \\
&  \textrm{with first non-zero eigenvalue denoted $\ai$}.
\end{aligned}
\right.
\end{equation*}
Now we  state a  result about the existence, the uniqueness and
convergence to equilibrium or the solution of (\ref{VPFPAintro}).
We call here solution on $[0,T[$ a function  $f \in \cc([0,T[,
\Binf^2)$ such that $\norm{E}_{L^\infty([0,T[\times \R^d)}<\infty$
and
\begin{equation*}
f(t) = e^{-tK}f_0 + \int_0^t e^{-(t-s)K}  E(s) \D_vf(s) ds.
\end{equation*}
We call it a global solution if it is a solution for all $T>0$.

\begin{thm} \label{maindecay}
 Suppose that conditions \Hun is satisfied. Then equation
 {\rm(\ref{VPFPAintro})}
 has a unique global solution for a given $L^1$-normalized
  initial datum $f_0 \in \Binf^2$.

  Besides if in addition \Hdeuxbis is fulfilled, then there exist
 constants $\Ai$ and $C_\infty$ only depending on second
order derivatives of $V + \Vinf$ and $\gamma$ (and uniform in
$\kappa$ varying in a fixed compact set) such that for any $\kappa
\leq \ai/C_\infty$
\begin{equation*}
\norm{ f(t,\cdot) -  \mminf}_{\Binf^2} \leq 6 \norm{ f_0 -
\mminf}_{\Binf^2} e^{-\frac{\ai }{2\Ai}t}.
\end{equation*}
\end{thm}

\bigskip
 As a corollary of Theorem \ref{maindecay}  we also get the
following result concerning the decay of the relative entropy.

\begin{cor} \label{cordecay}
Consider the solution given by Theorem {\rm\ref{maindecay}}. Then
(with the notations of Theorem {\rm\ref{maindecay}} and in
particular
  for $\kappa \leq \ai/C_\infty$) we have
\begin{equation*}
0\leq H(f(t),\mminf) \defegal \iint f(t) \ln
\sep{\frac{f(t)}{\mminf}}dxdv \leq C'_\infty \norm{ f_0
}_{\Binf^2}\norm{ f_0 - \mminf}_{\Binf^2}  e^{-\frac{\ai
}{2\Ai}t}.
\end{equation*}
where $C'_\infty$ only depends on second order derivatives of $V +
\Vinf$ and $\gamma$ (and is uniform in $\kappa$ varying in a fixed
compact set).
\end{cor}

Considering the short time linear diffusion estimates for
hypoelliptic operators, we mention the  cases $V = 0$ known since
\cite{Hor67} (see also the computations in the case  $V= x^2$ in
\cite{Kag01}) where the Green function is explicit. Numerous
non-linear result already quoted use this fact. For generic
hypoelliptic operators, this was studied by many authors in the
selfadjoint case,  in the spirit of the study of  sum of squares
of vector fields theorem with underlying Lie group structure. We
refer to the  book \cite{CSV92} and references therein for this
subject and point out that it is linked with the subelliptic
estimates for semi-groups of operators. The author was unable to
find any general result concerning the non-selfadjoint case (type
II operators), and the estimates given in Proposition \ref{disp0}
in this article seem to be new. Concerning
 the general study of the semi-group of globally hypoelliptic
  operators we also
 mention the recent works
 \cite{HN04},\cite{HelN04},  \cite{EH03},  \cite{HSS05}
  and  \cite{BL04}

Concerning the long time behavior of Fokker-Planck type operators,
we mention \cite{Tal99} for the exponential decay with probalistic
method, and \cite{HN04}  for the explicit exponential decay using
hypoelliptic tools developed first in \cite{EPR99}, \cite{EH00}.
Recent results  about more general kinetic equations on the Torus
can be found in works by Guo (e.g \cite{Guo02}).  We quote
\cite{DV01} for the convergence in $t^{-N}$ for all $N$ with the
use of entropy-dissipating methods, and  \cite{HSS05} for explicit
exponential decay using hypoelliptic tools close to the ones in
this paper. Let us also mention the work  \cite{Kag01} where
invariant manifolds methods are used in the case without external
potential.

 The systematic use of crossed derivatives in order to
get short-time, long-time and derivative estimates follows from
\cite{HN04} and was also developed in \cite{Vil05}-\cite{MN05}.
These studies, concerning type II operators follow common ideas
and features sometimes called hypocoercivity.

For the Vlasov-Poisson-Fokker-Planck equation there is a huge
literature on the subject (e.g. \cite{Deg86}, \cite{Bou93},
\cite{OV90}, \cite{RW92}  and \cite{OS00}). Essentially when $d=3$
these results
 use the explicit Green function and  $L^p$ estimates
available in this case. The case of a general confining potential
was not studied and in fact $L^p$ diffusion estimates on the
semi-group seem to be hard to get in this case. This is the reason
why we only deal with a mollified equation in the last part of
this paper. For the trend to the equilibrium, we quote
\cite{CSV96}, \cite{BD95},  \cite{Sol97} and  \cite{Dol99}.

\bigskip
The plan of the article is the following. In Section 2, we give
some results of functional analysis to be used later, essentially
taken from \cite{HN04}-\cite{HelN04}.  In Section 3 we prove
Theorem \ref{disp0} about the short-time diffusion estimate for a
general Fokker-Planck operator $K$.  There is a similar gain as in
the explicit case when the Green function is known (see e.g.
\cite{Bou93}) and obtained through hypoelliptic techniques. It
will play a crucial role in the study of the mollified VPFP
equation and in particular close to the equilibrium, where the
potential in $V+\Vinf$ is \it not known\rm. In Section 4 We give a
new proof of the exponential time decay toward the equilibrium,
based on an abstract Hilbert lemma given in the first subsection
there.  In the last section we apply the linear results first in
the case when an additional external field is added. Then we deal
with the Mollified VPFP equation and  prove Theorem
\ref{maindecay} and Corollary \ref{cordecay}.

\tableofcontents

\bigskip\bigskip

\noindent \bf Acknowledgement: \rm The author wants to thank C.
Villani and F. Nier for useful discussions, and D. Serre for
pointing out a mistake in a first version of this paper

\Section{Functional analysis} \label{functional}

We work here with a  potential $V$ independent of time and
satisfying conditions (H1). We denote by $B^2$ the  space defined
in (\ref{defb2}), and recall that it is an Hilbert with respect to
the scalar product
\begin{equation} \label{duality}
\seq{f,g} = \iint f g \mm^{-1} dxdv = \iint \frac{f}{\mm}
\frac{g}{\mm} \mm dxdv,
\end{equation}
for adequate $f$ and $g$. Recall that the spaces $\cc_0^\infty$
and $\mm^{1/2} \ss$ are dense in the $B^2$

We now state some results about the linear Fokker-Planck operator,
say in $L^1$. First the Hamiltonian vector field of $v^2/2 + V(x)$
is denoted by
$$
X_{0}=v\D_{x}-\D_{x}V(x)\D_{v},
$$
and it is easy to check that it is formally skew-adjoint with
respect to the scalar product (\ref{duality}) since $X_0$ commutes
with the multiplication with $\mm$. We also introduce the
differential ((d,1)-matricial) operators
\begin{equation} \label{ab}
a =\gamma^{1/2} (\D_x + \D_xV(x)), \ \ \ b =\gamma^{1/2} (\D_v +
v).
\end{equation}
For the scalar product defined in (\ref{duality}), their formal
adjoint are the following $(1,d)$-matricial operators
\begin{equation} \label{ab*}
a^* =  -\gamma^{1/2} \D_x,  \ \ \ b^* =-\gamma^{1/2} \D_v.
\end{equation}
With these notations the Fokker-Planck operator and its adjoint
with respect to the scalar product (\ref{duality}) read
\begin{equation*} \label{defK}
 K = X_0+ b^*b, \ \ \ \ K^* =-X_0+ b^*b.
\end{equation*}
Recall also  the definition of  Witten Laplacian (on $0$-forms) in
velocity and spatial coordinates
\begin{equation*}
\Lambda^2 = a^*a + b^*b,
 \end{equation*}
which is the naturally associated formally self-adjoint operator.
All these operators are  linked thanks to the following remarkable
algebraic properties:
$$
a = [b,X_0], \ \ \ b= -\Hess V [a,X_0].
$$
We want  to study the  linear Cauchy problem,
\begin{equation*} \label{lcp}
 \D_t f + K f =0, f|_{t=0} = f_0
 \end{equation*}
in  $B^2$. We first quote some  results from \cite{HelN04} and
\cite{HN04}.

\begin{prop}[\cite{HelN04}-\cite{HN04}] \label{cauchy} Operators $K$
and  $K^*$ defined as the closure of (\ref{defK})  with domain
$\cc^\infty_0$ are maximal accretive. They define semi-groups of
contraction and positivity preserving denoted $e^{-tK}$ (resp.
$e^{-tK^*}$) .
\end{prop}
We shall also need the following chain of Sobolev spaces based on
$B^2$. In the spirit of \cite{HN04} we denote
$$
\Lambda_a^2 = 1 + a^*a, \ \ \ \Lambda_b^2 = 1 + b^*b,
$$
where $a$ and $b$ were defined in (\ref{ab}-\ref{ab*}).
 Operators $\Lambda^2$,
$\Lambda_a^2$ and $\Lambda_b^2$  are maximal accretive with
$\mm^{1/2} \ss$ as a core and we denote by the same letter their
closure in $B^2$ (see the reference already quoted for instance).
In this sense $b^*b$ is an harmonic oscillator and $a^*a$ is the
Witten Laplacian associated to $V$. We introduce the natural chain
of Sobolev space for $k,l \in \R$
\begin{equation} \label{defhlk}
\H^{l,k} = \set{ f \in \mm^{1/2}\ss' \ \textrm{ s.t. }
\Lambda_a^k\Lambda_b^l f \in B^2 },
\end{equation}
for which $l\leq l' $ and $k\leq k'$ imply $\H^{l,k}
\hookrightarrow \H^{l',k'}$ and $\H^{0,0}=B^2$. We first write a
result which proof is essentially contained in
\cite{HN04}-\cite{HelN04} about the parabolic (smoothing)
properties of operator $K$ in $B^2$.

\begin{prop}[\cite{HN04}-\cite{HelN04}] \label{propregul}
For all $t > 0$ and $k\in\R$, $e^{-tK}$ maps $\mm^{1/2}\H^{-k,-k}$
to $\mm^{1/2}\H^{k,k}$. Besides for a given $k \geq 0$ there is
constants $C_{k,k}$ and $N_k$ such that for any initial data $f_0
\in \H^{k,k}$ we have
\begin{equation*} \label{regul}
\norm{f(t)}_{\H^{k,k}} \leq C_{k,k}\sep{t^{N_k} + t^{-N_k}}
\norm{f_0}_{\H^{-k,-k}}.
\end{equation*}
Besides the same result holds for  $K^*$.
\end{prop}

\remark \label{difference} The proof of this result is included in
\cite{HN04}. Let us just notice some differences. Here there is no
growing assumptions for $\D_x V$, anyway the definition of the
commutators and the pseudodifferential calculus are still valid
with based metric $g_0 = dx^2 + d\xi^2 + d\eta^2 + dv^2$ where
$(\xi,\eta)$ are the dual variables of $(x,v)$. On the contrary
some assertions about compactness (of the resolvent of  $K$,
$\Lambda^2$ ...) are not true anymore. It corresponds in
\cite[Appendix A]{HN04} to the case $n=1/2$. In particular
$\cap_{k\in\R} \H^{k,k} \neq \ss$ and $e^{-tK}$ does not anymore
send $\mm^{1/2}\ss'$ to $\mm^{1/2}\ss$. Anyway the proofs of the
other result there remain true under the hypothesis \Hun here.

\Section{Short time behavior}

The purpose of the following Section is to prove Theorem
\ref{disp0} about the short time behavior of the semi-group
associated to $K$. In particular we want to ameliorate the
estimate for small $t$ in (\ref{regul}), at least in the case $k
=1$, and with explicit bounds. It is based on the construction of
a particular Lyapounov functional $A(t)$ taking into account the
hypoelliptic properties of $K$. In fact Theorem \ref{disp0} is
included in the following Proposition:

\begin{prop} \label{disp}
There exists a constant $C_2$  such that
 for all $ t > 0 $, we have the following: \\
 i)  $e^{-tK} b^* \textrm{ is bounded by } C_2 (1+t^{-1/2})
$ and\\
ii) $ e^{-tK} a^* \textrm{ is bounded by } C_2 (1+t^{-3/2}),
$ \\
as bounded operators on $B^2$. Here  $C_2$ depends only
 on $\norm{V''}_{L^\infty}$ (and $\gamma$). Besides we have
 the same bounds as in i) for the  operators $b^\natural
 e^{-tK^\sharp}$ and $e^{-tK^\sharp} b^\natural$  and as in ii) for
 operators $a^\natural
 e^{-tK^\sharp}$ and $e^{-tK^\sharp} a^\natural$, where $\sharp$
 and $\natural$ are either nothing or $*$.
\end{prop}

\preuve We shall in a moment prove the results for $ae^{-tK}$ and
$be^{-tK}$. Taking these bounds for given we note that they imply
similar bounds for their  adjoints $e^{-tK^*}b^*$ and
$e^{-tK^*}a^*$ since $ {B^{2}}'=B^2$. The proof is
 exactly the same for $be^{-tK^*}$ and $a e^{-tK^*}$
  since the sign
 in front of $X_0$ has essentially no importance in the proof.
 Taking the adjoints
 again give the result for $e^{-tK}b^*$ and
$e^{-tK}a^*$.

Now for the bound on $b^*e^{-tK}$ we simply have to note that for
$f_0 \in B^2$ given  and $f(t) = e^{-tK}f_0$  we have
$$
\norm{b^*f(t)}^2 = ( bb^* f(t), f(t)) = ( b^*b f(t), f(t)) +
d\norm{f(t)}^2 = \norm{bf(t)}^2 + d \norm{f(t)}^2
$$
and we get the result. For $a^*e^{-tK}$ we similarly write
$$
\norm{a^*f(t)}^2 = ( aa^* f(t), f(t)) = ( a^*a f(t), f(t)) +
(\Delta \V f(t),f(t)) \leq   \norm{af(t)}^2 + C \norm{f(t)}^2
$$
since $\V$ is with second derivatives bounded. This gives the
result. For the other terms, we repeat the procedure followed in
preceding paragraph and the proof of the last assertion in
Proposition \ref{disp} is complete.

\bigskip Let us come back now to the bounds on $ae^{-tK}$ and
$be^{-tK}$. We  note that the operators are well defined since
$e^{-tK}$ is defined from $\mm^{1/2}\H^{-k,-k}$ to
$\mm^{1/2}\H^{k,k}$. For the proof we shall need a series of
results.

\bigskip
 The first thing we do it to change the function $f$ by the
standard conjugation tool: We pose
\begin{equation} \label{flatu}
u = f/\mm^{1/2}, \ \ \ u_0 = f_0 / \mm^{1/2}.
\end{equation}
After this conjugation, operator $K$ is replaced by the following
\begin{equation} \label{flatK}
 K = v\D_x - \D_xV \D_v + \gamma( -\D_v + v/2)(\D_v + v/2).
\end{equation}
acting on $u$, and defined in the flat space $L^2$.  We recall
 also that the Witten Laplacian in both variables $v$ and $x$
 reads after conjugating
\begin{equation} \label{flatlambda}
 \Lambda^2 = \gamma(-\D_x + \D_x V/2).(\D_x + \D_x V/2)
 + \gamma(- \D_v +v/2)(\D_v + v/2),
\end{equation}
and that operators $a$, $b$ and their adjoints are now
\begin{equation} \label{flatab}
\begin{array}{lll}
& a= \gamma^{1/2}(\D_x + \D_x V/2), \ \ \ & a^* =
\gamma^{1/2}(-\D_x
+ \D_x V/2), \\
\ \ \textrm{ and } \ \ & b= \gamma^{1/2}(\D_v +v/2), \ \ \ & b^* =
\gamma^{1/2}(- \D_v +v/2).
\end{array}
\end{equation}
Operator $X_0 = v\D_x - \D_x V \D_v $ is unchanged and we have
again
\begin{equation} \label{flat3}
K = X_0 + b^*b, \ \ \ \Lambda^2 = a^*a + b^*b.
\end{equation}
For all this conjugated operators we keep the same notations as
for the unconjugated ones. No confusion is possible since they act
in $L^2$ on the conjugated function $u$ instead of $f$. The norm
is the standard one associated with the $L^2$ space.

\bigskip We work now in
 the $L^2$ setting we just defined.
 We recall that for any $l\in \R$ and   $u_0 \in \Lambda^{l} L^2$,
  Proposition \ref{propregul} implies that for all $t>0$  and $k\in
  \R$,
$u(t)\in \Lambda^kL^2$. We now choose $u_0$ such that $\Lambda^2
u_0 \in L^2$ and we pose for $t \geq 0$,
$$
A(t) = t^3\norm{a u }^2 + E t^2 \Re (au,bu) + D t \norm{bu}^2 + C
\norm{u}^2.
$$
 $A$ is a   $\cc^0(\R^+, \R) \cap \cc^1(\R^{+*}, \R)$ function,
 and we can compute its time derivative for $t>0$.

 \ititem{Derivative of $\norm{u}^2$:}
 We have
 \begin{equation} \label{ordre0}
 \D_t \norm{u}^2 = -2\Re(Ku,u) = -2 \norm{bu}^2.
 \end{equation}

\ititem{Derivative of $t\norm{bu}^2$:} We write
\begin{equation*} \label{ordre1temp1}
 \D_t t\norm{bu}^2 =  \norm{bu}^2 + t \D_t(b^*bu,u).
 \end{equation*}
Let us compute separately  the second derivative. We have
\begin{equation*}\begin{split}
\D_t(b^*bu,u) & = -\Re (b^*b K u ,u) -\Re(b^*bu,Ku) \\
 & = -2 \norm{b^*bu}^2 -\Re (b^*b X_0 u ,u) + \Re(X_0 b^*bu,u) \\
 &  = -2 \norm{b^*bu}^2 - \Re \sep{ \adf{b^*b, X_0}u,u}.
 \end{split}
 \end{equation*}
 using that $a = [b,X_0]$, we get
\begin{equation*}\begin{split}
 \D_t(b^*bu,u)& = -2 \norm{b^*bu}^2 - \Re (b^*au,u)
  - \Re(a^*b u,u) \\
 & = -2 \norm{b^*bu}^2 - 2\Re(au,bu).
 \end{split}
 \end{equation*}
 As a consequence we can write that
\begin{equation} \label{ordre1}
\D_t( t\norm{bu}^2) = \norm{bu}^2 -2t \norm{b^*b u }^2 -2 \Re
t(au,bu).
\end{equation}

\ititem{Derivative of $t^2 \Re(au,bu)$:} We write
\begin{equation} \label{ordre2temp1}
 \D_t t^2\Re\sep{au,bu} =  2t\Re(au,bu) + t^2 \D_t\Re(au,bu).
 \end{equation}
Let us compute again separately  the second derivative :
\begin{equation*}\begin{split}
\D_t\Re(au,bu) & = -\Re(aKu,bu) -\Re(au, bKu) \\
&= -\Re(ab^*bu,bu) -\Re(au, bb^*bu) -\Re(a X_0u, bu)
-\Re(au,bX_0u).
\end{split}
 \end{equation*}
We can commute the field $X_0$ in the last two terms and we get
\begin{equation*}\begin{split}
\D_t\Re(au,bu)  = &  -\Re(ab^*bu,bu) -\Re(au, bb^*bu)
 -\Re([a, X_0]u,
bu)
-\Re(au,[b,X_0]u) \\
& + \underbrace{\Re( X_0 au, bu) + \Re(au,X_0 bu)}_{=0 \text{
since} X_0 \text{ is skewadjoint }}.
\end{split}
 \end{equation*}
Now use the facts that $ [b, X_0]= a$ and $ -\Hess V b = [a,X_0]$.
This yields
\begin{equation*}\begin{split}
\D_t\Re(au,bu)  = &  -\Re(bau,bbu) -\Re(b^*au, b^*bu)
 +(\Hess V bu,bu) - \norm{au}^2
\end{split}
 \end{equation*}
and using  (\ref{ordre2temp1}) we get
\begin{equation*} 
\begin{split}
\D_t \sep{ t^2\Re \sep{au,bu} } = & 2t \Re(au,bu) - t^2\norm{au}^2
+t^2(\Hess V bu,bu) \\
&  -t^2\Re(bau,bbu) -t^2\Re(b^*au, b^*bu).
\end{split}
 \end{equation*}
Using eventually the fact that $b^*b = bb^* -\gamma d$ yields
\begin{equation} \label{ordre2}
\begin{split}
\D_t \sep{ t^2\Re \sep{au,bu} } = & 2t \Re(au,bu) - t^2\norm{au}^2
+t^2(\Hess V bu,bu) \\
&  -2t^2\Re(bau,bbu) -t^2 \gamma d\Re(au, bu).
\end{split}
 \end{equation}

\ititem{Derivative of $t^3 \norm{au}^2$:} We write
\begin{equation} \label{ordre3temp1}
\D_t t^3 \norm{au}^2 = 3 t^2 \norm{au}^2 + t^3 \D_t \norm{au}^2.
\end{equation}
We study separately the second term:
\begin{equation*}\begin{split}
\D_t\Re(au,au) & = -\Re(aKu,au) -\Re(au, aKu) \\
&= -\Re(ab^*bu,au) -\Re(au, ab^*bu) -\Re(a X_0u, au)
-\Re(au,aX_0u).
\end{split}
 \end{equation*}
We shall again commute the field $X_0$ in the last terms
\begin{equation*}\begin{split}
\D_t\Re(au,au)  = &  -2\norm{bau}^2 -2\Re(au, aX_0u) \\
= & -2\norm{bau}^2 -\Re([a, X_0]u,au) + \underbrace{\Re( X_0 au,
au)}_{=0 \text{
since} X_0 \text{ is skewadjoint }} \\
\end{split}
 \end{equation*}
Now since $ -\Hess V b = [a,X_0]$. We get
\begin{equation*}\begin{split}
\D_t\Re(au,au)  = &  -2 \norm{bau}^2
 +2\Re(\Hess V bu,au).
\end{split}
 \end{equation*}
From (\ref{ordre3temp1}) we can therefore write
\begin{equation} \label{ordre3}
\begin{split}
\D_t t^3 \norm{au}^2 = 3 t^2 \norm{au}^2 -2 t^3 \norm{bau}^2 +2t^3
\Re(\Hess V bu,au).
\end{split}
 \end{equation}

\ititem{Derivative of $A$:} We put together the results of
(\ref{ordre0}-\ref{ordre1}-\ref{ordre2}-\ref{ordre3}) and we get
the following formula for the derivative of $A$, where we have put
the similar terms on the same lines :
\begin{equation*}
\begin{array}{rlr}
&\D_t A(t) = & \\
 & -2C \norm{bu}^2 -2tD\norm{b^*bu} -t^2E
\norm{au}^2 -
2t^3 \norm{bau}^2 & \entoure{1}\\
& +D\norm{bu}^2  + t^2 E\Re(\Hess V bu,bu) & \entoure{2} \\
& + 2 tD\Re(au,bu) + 2tE\Re(au,bu) + 2 t E \Re(\Hess V bu,au) -t^2
D\gamma d\Re(au, bu)
& \entoure{3} \\
& + 3 t^2 \norm{au}^2 & \entoure{4} \\
& -2E t^2 \Re(bau,bbu). & \entoure{5}
\end{array}
\end{equation*}
We bound now each terms on the lines \entoure{2}-\entoure{5} by
terms appearing in \entoure{1}. We suppose that $t\in ]0,1]$. Now
since the Hessian of $V$ is bounded  by a constant, say $C_V$, we
have
\begin{equation} \label{bound2}
\entoure{2} \leq (D + E C_V) \norm{bu}^2   \ll 2C \norm{bu}^2
\text{if } D,E \ll C.
\end{equation}
For the term \entoure{3}, we write for $\eta >0$,
\begin{equation*}
\begin{split}
\entoure{3} &\leq (2D+2E + 2EC_V + D\gamma n)
 t\norm{au}\norm{bu} \\
&  \leq \eta t^2 \norm{au}^2 + \frac{Cte(D,E,\gamma d)}{\eta}
\norm{bu}^2.
 \end{split}
 \end{equation*}
 We therefore get that for a given $E$, we have to choose
  $\eta$ sufficiently small
 and then  $C$ big enough to get
\begin{equation} \label{bound3}
\begin{split}
\entoure{3} &  \ll E t^2 \norm{au}^2 + 2C \norm{bu}^2.
 \end{split}
 \end{equation}
Now we treat the term \entoure{4} : this is easy since we only
need to take $E \gg 3$ in order to get
\begin{equation} \label{bound4}
\entoure{4} \ll -3t^2\norm{au}^2.
\end{equation}
For the last term \entoure{5} we write :
\begin{equation*}
\begin{split}
\entoure{5} =
-2Et^2\Re(bau,bb u) & \leq 2E t^2\norm{ba u}\norm{bbu} \\
& \leq E( \eta' t^3 \norm{bau}^2 + \frac{t}{\eta'} \norm{bb u}) \\
& \leq  E \eta' t^3 \norm{bau}^2 + \frac{Et}{\eta'} \norm{b^*b u},
\end{split}
\end{equation*}
where in the last estimate we use the fact that for $w$ s.t.
$\Lambda_b w \in L^2$, we have $\norm{bw} \leq \norm{b^*w}$. Now
for a given $E$ we have to choose first $\eta'$ small enough, and
then $D$ sufficiently large to write
\begin{equation} \label{bound5}
\entoure{5}   \ll 2tD\norm{b^*bu}^2 + 2t^3\norm{bau}^2.
\end{equation}

\ititem{Synthesis} We checked that each line
\entoure{2}-\entoure{5} can separately be bounded by a term
appearing in \entoure{1}. In order to get the fact that $t\mapsto
A(t)$ is decreasing, we choose the constants as follows : first
 $E$ so that (\ref{bound4}), and then  $\eta$,
$\eta'$, $C$ and $D$ such that (\ref{bound5}) and (\ref{bound3}).
Eventually increasing $C$ so that  (\ref{bound2}) holds yields the
result. In particular since $A(t)$ is right-continuous in $0$ we
get that for all $t \in [0,1]$
\begin{equation*}
A(t) = t^3\norm{a u }^2 + E t^2 \Re (au,bu) + D t \norm{bu}^2 + C
\norm{u}^2 \leq C\norm{u_0}^2.
\end{equation*}
In particular we have for $t\in [0,1]$,
\begin{equation*} \label{short}
\norm{au(t)} \leq C^{1/2}t^{-3/2} \norm{u_0}, \ \ \ \ \norm{bu(t)}
\leq  (C/D)^{1/2} t^{-1/2} \norm{u_0}.
\end{equation*}
This is the short time  estimate ($t \in ]0,1]$) in Proposition
\ref{disp} for $\Lambda^2u_0 \in L^2$. For $t\geq 1$ we simply
write that
$$
\norm{be^{-tK}u_0} = \norm{ be^{-K/2} e^{-(t-1/2)K} u_0} \leq C_2
\norm{e^{-(t-1/2)K} u_0} \leq C_2 \norm{u_0},
$$
where we used first the short time estimate (with $t=1/2$) and
then the fact that $K$ is maximal accretive. The result for $u_0
\in L^2$ follows then by density. \fin

\Section{Exponential time decay} \label{secdecay}

\Subsection{An abstract Hilbert result} \label{explindecay}

Let us first state a general lemma about semi-group of operators.
Let $\kk$ be the infinitesimal generator of a semi-group of
contraction on a Hilbert space $H$ (in particular
$\overline{D(\kk)}= H$). We want to extend the following basic
result :
\begin{quote}
If $\exists \alpha >0$ such that $\alpha \norm{\phi}^2 \leq
\Re(\kk \phi,\phi)$ for all $\phi\in D(\kk)$ \\
then $\forall \phi_0 \in H, \ t\geq 0$ we have $\norm{e^{-t\kk}
\phi_0} \leq e^{-\alpha t} \norm{\phi_0}$,
\end{quote}
Of course the converse is true  applying the Lummer Phillips
theorem to the operator $\kk-\alpha \Id$ (see for example
\cite{Paz83}). We want to extend the right sense.

\begin{lem} \label{expdecayabs}
Let $\kk$ be the infinitesimal generator of a semigroup of
contraction on a Hilbert space $H$ and suppose that there exist a
constant  $\alpha >0$ and a bounded operator  $\lll$ with norm
bounded by  $C \geq 1$
 such that
\begin{equation} \label{abstr}
\forall \phi  \in D(\kk), \ \ \ \alpha\norm{\phi}^2 \leq \Re(\kk
\phi,\phi) + \Re(\kk\phi, (\lll+\lll^*)\phi)
\end{equation}
then for all $\phi_0 \in H$ and $t \geq 0$ we have
\begin{equation*} \label{abslem}
\norm{ e^{-t\kk} u_0} \leq 3 e^{-\frac{\alpha t }{3C}}
\norm{\phi_0}
\end{equation*}
\end{lem}

\preuve We write for $\phi_0 \in D(\kk)$, $\phi(t) = e^{-t\kk}
\phi_0 \in D(\kk)$. Using (\ref{abstr}) and since $\kk$ is
accretive we get
\begin{equation*}
 3 C \frac{\alpha}{3C} \norm{\phi}^2 \leq 4 C \Re(\kk\phi,\phi) +
 \Re(\kk\phi,(\lll+\lll^*)\phi)
 \end{equation*}
 and  since $|\Re(\lll\phi,\phi)| \leq C\norm{\phi}^2$ we have
\begin{equation*}
  \frac{\alpha}{3C}\sep{ 2C \norm{\phi}^2 +
  \Re(\lll\phi,\phi)} \leq  4 C \Re(\kk \phi,\phi) +
 \Re(\kk \phi,(\lll+\lll^*)\phi).
 \end{equation*}
 Now $\D_t \norm{\phi}^2 = -2\Re(\kk\phi,\phi)$ and $\D_t \Re(\lll\phi,\phi) =
 -\Re(\kk\phi, (\lll+\lll^*)\phi)$ therefore
\begin{equation*}
  \frac{\alpha}{3C}\sep{ 2C \norm{\phi}^2 + \Re(\lll\phi,\phi)} +
  \frac{\D}{\D_t} \sep{ 2 C \norm{\phi}^2+
 \Re(\lll\phi,\phi)} \leq 0.
 \end{equation*}
 Integrating between $0$ and $t$ gives
 $$
 2 C \norm{\phi}^2+
 \Re(\lll\phi,\phi) \leq e^{-\frac{\alpha t}{3C}}
 \sep{2 C \norm{\phi_0}^2+\Re(\lll\phi_0,\phi_0)}.
 $$
 Using twice the fact that $\norm{\lll}$ is bounded by $C$ we get
 $$
 C \norm{\phi}^2 \leq 3C e^{-\frac{\alpha t}{3C}}\norm{\phi_0}^2.
 $$
 This gives (\ref{abslem}) since $\norm{\phi} \leq \norm{\phi_0}$  and
 the proof is complete. \fin

\subsection{The case of the Fokker-Planck operator}
\label{lindecay}

We want now to apply the preceding abstract result to the linear
Fokker-Planck equation and prove Theorem \ref{mainn}. A part of
the proof is very close to the one given in \cite{HN04} and in
particular uses Kohn's type arguments about hypoellipticity
developed there. We work from now on with a
 potential $V$ satisfying both conditions of type (H1) and
(H2).

We shall work in the following in the orthogonal in the $B^2$
sense of the Maxwellian. For the following we call
$$
\Bperp^2 = \mm^{\perp} \cap B^2 = \set{f \in B^2 \text{ s.t. }
\int f dxdv = 0},
$$
 endowed with the norm of $B^2$,
  where $\perp$ stands for the orthogonal with respect to
the scalar product in $B^2$ (recall that $B^2{}'$ was identified
with $B^2$ according to the measure $\mm^{-1}dxdv$ in
(\ref{duality})).

We note that $\Bperp^2$ is stable for $K $. Indeed for all $f \in
\Bperp^2 $, we have
$$
<K f,\mm> = \int K f \mm \mm^{-1} dx = \int f (K^* \mm) \mm^{-1} =
0.
$$
Since $K$ restricted to $B^2$ generates a semi-group of
contraction, we have the same property in $\Bperp^2$. Anyway
restricted to $\Bperp^2$ the semi-group has a much better property
at infinity:

\begin{prop} \label{expdecay}
Suppose  $V$ satisfy hypothesis \Hun and \Hdeux. Then there exists
a constant $A$  depending only on
 $\norm{V''}_{L^\infty}$ such that for all $t\geq 0$
$$
\norm{e^{-tK}}_{\Bperp^2} \leq 3 e^{-\alpha t / A}
$$
where $\alpha$ was defined in \Hdeux. The same bound occurs for
$K^*$.
\end{prop}

\preuve Again we work in $L^2$ after conjugation  by the square
root of the Maxwellian, and the new unknown function is $u$. We
therefore keep for the proof of the Proposition the notations $K$,
$a$, $b$, $X_0$ and $\Lambda^2$ introduced in
(\ref{flatu}--\ref{flat3}).
 Note that after the conjugation, the space $B^2$ becomes the orthogonal
 of the square root of the Maxwellian and we denote it by
 $H= (\mm^{1/2})^\perp$. We note that $H$ is stable
for both $K$ and $\Lambda^2$ and that $\mm^{1/2}$ is a single
eigenfunction for $K$ and $\Lambda^2$ with eigenvalue $0$. We also
introduce the
 following operator on (the flat space ) $L^2$:
$$
 \Lambda^2_\delta
\defegal \delta^2 + a^*a + b^*b.
$$
where $\delta^2 \leq \gamma$ is to be fixed later.

Let's take $u \in \ss$. We first quote the result of Proposition
2.5 (case $\eps = 0$ there) in \cite{HN04}, which is true under
our assumptions on $V$.
\begin{equation}
\begin{split}
\label{eq.lem2} \norm{ u}^2 \leq \ & \Re \sep{K u,\sep{L + L^*} u
}
 -  2\Re \sep{ b^*b u ,L u }  -\Re(\aaa^*bu,u) \\
& + (1+\gamma) \delta^{-1} \norm{b u}\norm{u}  + \delta^2
(\Lambda^{-2}_\delta u,u),
\end{split}
\end{equation}
 where $L = \Lambda^{ -2}_\delta a^* b$ and
 $\aaa^*=\adf{ \Lambda^{-2}_\delta a^*, X_0}$.
From Proposition 5.4 in \cite{HN04} we have an explicit bound for
the bounded operator $\aaa^*$, and this is also easy to get bounds
for $L$ and $a\Lambda^{-2}_\delta b^*$,
\begin{equation*} \label{CVVV}
\norm{\aaa^*} \leq C_V \delta^{-1}, \ \ \ \norm{L} \leq
\sqrt{2d\gamma} \delta^{-1}, \ \ \ \norm{a\Lambda^{-2}_\delta b^*}
\leq 1.
\end{equation*}
(For the second one we simply observe that if $a_j$, $b_j$ denote
the components of $a$ and $b$ we have $ \norm{b^*au} \leq \sum_j
\norm{b_j^*a_j u}$ and $\norm{b_j^*a_j u}^2 = (b_jb_j^* u ,
a_j^*a_j u ) = (b_j^*b_ju,a_j^*a_j u ) + \gamma(  a_j^*a_j u,u )$
and since $\delta^2 \leq \Lambda^2_\delta$ we get  $\norm{b_j^*a_j
u}^2 \leq 2\gamma\delta^{-2} \norm{\Lambda^2_\delta u}^2$
therefore $ \norm{b^*au}^2 \leq 2d\gamma\delta^{-2}
\norm{\Lambda^2_\delta u}^2$ and $ \norm{b^*a \Lambda^{-2}u}^2
\leq 2d\gamma\delta^{-2} \norm{ u}^2$. Taking the adjoint and the
square root gives the result.)

We can then write from (\ref{eq.lem2}) that
\begin{equation*}
\begin{split}
\norm{ u}^2 \leq \ & \Re \sep{K u,\sep{L + L^*} u }
 +  2 | \sep{ a\Lambda^{-2}_\delta b^*b u ,b u }|  + |\Re(\aaa^*bu,u)| \\
& + (1+\gamma) \delta^{-1} \norm{b u}\norm{u}  + \delta^2
(\Lambda^{-2}_\delta u,u) \\
\leq & \Re \sep{K u,\sep{L + L^*} u } + 2 \norm{bu}^2 + C_V
\delta^{-1}\norm{bu}\norm{u} \\
&  + (1+\gamma)\delta^{-1} \norm{b u}\norm{u}  + \delta^2
(\Lambda^{-2}_\delta u,u).
\end{split}
\end{equation*}
Using first the inequality $|xy| \leq x^2 + 4^{-1}y^2$ and then
the fact that $\Re(Ku,u) = \norm{bu}^2$ we get
\begin{equation} \label{presquefin}
\begin{split}
\norm{ u}^2  \leq & \Re \sep{K u,\sep{L + L^*} u }  +
C_V'\delta^{-2}\norm{bu}^2
 + \frac{1}{4} \norm{u}^2 + \delta^2
(\Lambda^{-2}_\delta u,u) \\
\leq & \Re \sep{K u,\sep{L + L^*} u }  + C_V'\delta^{-2}\Re(Ku,u)
 + \frac{1}{4} \norm{u}^2 + \delta^2
(\Lambda^{-2}_\delta u,u).
\end{split}
\end{equation}
Now we suppose that $u \in H$ therefore $\delta^2
(\Lambda^{-2}_\delta u,u) \leq
\frac{\delta^2}{\alpha+\delta^2}\norm{u}^2$ and we choose
$\delta^2 = \alpha$ (which is lower than $\gamma$ because of the
harmonic part of $\Lambda_0^2$). This gives
$$
\delta^2(\Lambda^{-2}_\delta u,u) \leq \frac{1}{2} \norm{u}^2
$$
and putting this in (\ref{presquefin}) we get
\begin{equation*} 
\begin{split}
\frac{1}{4}\norm{ u}^2  \leq& \Re \sep{K u,\sep{L + L^*} u }  +
C_V'\delta^{-2}\Re(Ku,u).
\end{split}
\end{equation*}
As a consequence
\begin{equation*} \label{presquefin2}
\begin{split}
\frac{\alpha}{4C_V'}\norm{ u}^2  \leq & \Re \sep{K
u,\sep{\widetilde{L} + \widetilde{L}^*} u } + \Re(Ku,u)
\end{split}
\end{equation*}
where $\widetilde{L} = \frac{\delta^2}{C_V'}L$ satisfies
$\norm{\widetilde{L}} \leq \sqrt{2d\gamma} \delta^{-1} \delta^2
/C_V' \leq 1$ since $ \delta \sqrt{2d \gamma} \leq \sqrt{2d}\gamma
\leq C_V'$. The result of the lemma is then a direct consequence
 of Lemma \ref{expdecayabs}. Taking $A=12 C_V'$ completes the
 proof of the Proposition. \fin

\remark We can point out that the gain with respect to the
estimate in \cite{HN04} is first  that the constant in front of
the exponential is universal ($=3$) and in particular does not
depend on $V$ or $\alpha$. It makes sense in the study of the
Vlasov-Poisson-Fokker-Planck system with small data in Section
\ref{secdecay} since this constant has to be compared with the
size of the initial data. The second remark is that no assumption
about the increasing of $\D_x V$ is made, and we can understand
this fact by saying that the existence of a   spectral gap for the
Witten Laplacian implies  a (generalized) spectral gap for the
Fokker-Planck operator, without assumptions on the remaining part
 of the spectrum (implied for example by the compacity
 of the resolvent). \fin

 \remark Note to the end  that in the preceding study
is also valid even for $V''$ not in $L^\infty$ since the only real
needed assumption is that the constant $C_V$ in (\ref{CVVV}) is
not infinite (see \cite{HN04} for its expression). Anyway in this
case one has to be careful when defining  the commutators, and
some additional assumptions on $V$ may be needed. \fin

Putting together Propositions \ref{expdecay} and \ref{disp} we can
complete the proof of Theorem \ref{mainn}. In fact it is included
in the following Proposition

\begin{prop} \label{disp2}
Suppose  $V$ satisfy hypothesis \Hun and \Hdeux.
Then there exists  constants $C$ and $A$  such that for all $t>0$, \\
 i)\ $be^{-tK}$, $e^{-tK} b^*$, $be^{-tK^*}$ and $e^{-tK^*} b^*$
 \textrm{ are bounded by } $C (1+t^{-1/2})
e^{-\alpha t/A}$ \\
ii)\ $ae^{-tK}$, $e^{-tK} a^*$, $ae^{-tK^*}$ and $e^{-tK^*} a^*$
\textrm{ are bounded by } $C (1+t^{-3/2})
e^{-\alpha t/A}$ \\
as bounded operators in $B^2$, where $C$ and $A$ depends only on
$\norm{V''}_{L^\infty}$.
\end{prop}

\preuve[of Proposition \ref{disp2} and Theorem \ref{mainn}] We
simply use the following fact : For a given $f_0 \in \mm^{1/2} \ss
\subset B^2$, we have
$$
(b^*f_0 ,\mm)_{B^2}= (f_0,b\mm)_{B^2} =  0.
$$
i.e. $ b^*f_0$ is orthogonal to the square root of the Maxwellian.
 Of course it is also the case for $ e^{-K}b^*f_0$ since
$$
(e^{-K}b^*f_0 ,\mm)_{B^2}= (b^*f_0,e^{-K^*}\mm)_{B^2} =
(b^*f_0,\mm)_{B^2} =  0.
$$
Now for $t\geq 1$ we can apply  Proposition \ref{expdecay} and we
get  following bound :
\begin{equation} \label{bonbon}
\begin{aligned}
\norm{e^{-tK}b^*f_0} & = \norm{e^{-(t-1)K}e^{-K}b^*f_0} \leq C
e^{-\alpha (t-1)/A}\norm{ e^{-K}b^*f_0} \\
& \leq 3e^{\gamma/A} e^{-\alpha t}\norm{ e^{-K}b^*f_0},
\end{aligned}
\end{equation}
since  $\alpha \leq \gamma$ because of the harmonic part of
$\Lambda^2$. Now from Proposition \ref{disp} applied with $t=1$ we
get
$$
\norm{e^{-K} b^* f_0}_{B^2} \leq C'(1+t^{-1/2})\norm{f_0}_{B^2}.
$$
This inequality together with (\ref{bonbon}) give the bound
$$
\norm{e^{-tK} b^* f_0}_{B^2} \leq C_2(1+t^{-1/2})e^{-\alpha
t/A}\norm{f_0}_{B^2}
$$
 for an initial data $f_0 \in \mm^{1/2} \ss$.
It can be clearly extended to $f_0 \in B^2$ by density. The proof
of the estimates about $e^{-tK} a^* $ can be done exactly in the
same way. The assertions concerning $K^*$ are immediate since the
structure of the operator is the same. Eventually the ones
concerning $ae^{-tK}$ and $be^{-tK}$ are immediate using the
adjoint of the preceding ones. The proof of Proposition
\ref{disp2} is complete. \fin

\remark In the particular case of $V = x^2$ or more generally for
quadratic external potentials, one can compute explicitly the
Green function of $e^{-tK}$ using the method of characteristics
(see e.g. \cite{Kag01}). Anyway if it gives after some work the
short time behavior, the exponential decay of $e^{-tK}b^*$ is not
clear on the formulas. In fact the short time decays in theorems
(\ref{disp}) can be viewed as  consequences of the Lie  group
structure
 of the operator (if one assimilates $b^*$ and $b$) whereas
 the long time behavior is deeply linked with the spectral
 properties of $K$.

\Section{Applications}
 \Subsection{Strong solutions for a given
interaction potential} \label{degond}

In this section we work again in a linear setting  and study the
following Fokker-Planck equation
\begin{equation}
  \label{FPE}
  \left\{
  \begin{aligned}
   &  \D_{t}f+v.\D_{x}f -(E+ \D_x V).\D_{v}f
      -\gamma \D_{v}.\left(\D_{v}+v\right)f= U \\
   &   f|_{t=0} = f_0,
 \end{aligned}
 \right.
\end{equation}
where $E(t,x)$ is a given time-dependant potential satisfying $ E
\in L^\infty([0,T[ \times \R^d) $ and $V$ is again a potential
satisfying hypothesis (H1). We shall prove  existence and
uniqueness in the space $\H^{l,k}$ based on $B^2$ and defined in
(\ref{defhlk}). In  the following Propositions we will assume the
following
\begin{equation*}
\Htrois \ \ \ \ \left\{
\begin{aligned}
& V \textrm{ satisfies  (H1)}, \\
& E \in L^\infty([0,T[ \times \R^d), \\
& f_0 \in B^2(\R^{2n}), \\
& U \in L^2 \sep{ [0,T[ , \H^{0,-1}(\R^d_x\times \R^d_v)}.
\end{aligned}
\right.
\end{equation*}

The first Proposition gives  existence and uniqueness of a unique
mild solution of the system (\ref{FPE}).

\begin{prop} \label{existmild}
Under hypotheses   \Htrois  there exists a unique mild solution
solution of {\rm(\ref{FPE})}, where by definition a mild solution
is a solution $ f \in \cc([0,T[, B^2)$ satisfying
$$
f(t) = e^{-tK}f_0 + \int_0^t e^{-(t-s)K}  E(s)\D_v f(s) + \int_0^t
e^{-(t-s)K}  U(s) ds.
$$
\end{prop}

\preuve This is obtained via a  standard fixed point theorem
 in $L^\infty([0,T[, B^2)$. We only sketch the proof.
 To simplify the notations we suppose $\gamma=1$
  which implies $\D_v= -b^*$.
 Let $F$ be the following operator from
$L^\infty([0,T[, B^2)$ into itself given by
\begin{equation} \label{fpE}
\begin{split}
F(f) & = e^{-tK}f_0+ \int_0^t e^{-(t-s)K}  E(s) \D_v f(s) +
\int_0^t e^{-(t-s)K}
U(s) \\
& =e^{-tK}f_0- \int_0^t e^{-(t-s)K} b^* E(s) f(s) + \int_0^t
e^{-(t-s)K}(1+b^*b) \Lambda_b^{-1} \Lambda_b^{-1} U(s).
\end{split}
\end{equation}
According to the diffusion estimates given in Proposition
\ref{disp} and using the fact the $b\Lambda_b^{-1}$ is  bounded by
$1$ as an operator  in $B^2$, we get that for all $0<t<T$
$$
\norm{e^{-tK}(1+b^*b) \Lambda_b^{-1} }_{B^2} \leq Ct^{-1/2}.
$$
We therefore get  for $f\in B^2$, $F(f) \in B^2$ and  for
 all $t>0$,
$$
\norm{F(f)}_{L^{\infty}([0,t[ ,B^2)} \leq C t^{1/2} \norm{f
}_{L^{\infty}([0,t[ ,B^2)}.
$$
Using a standard fixed point theorem we get that $f$ is the unique
limit of the following iteration scheme
\begin{equation*} \label{sequencefn}
 \D_{t}f^{n+1} +K f^{n+1} +E b^* f^n
      = U, \ \ \ \ f^0 = f_0,
 \end{equation*}
and the continuity is clear from formula (\ref{fpE}). \fin

\subsection{A mollified Vlasov-Poisson-Fokker-Planck system}
 In this section we study the following non-linear problem, to be
 understood as a  modified Vlasov-Poisson-Fokker-Planck system, where the
 non-linear coupling is mollified:

\begin{equation}
  \label{VPFPA}
  \left\{
  \begin{aligned}
   &  \D_{t}f+v.\D_{x}f -(E+ \D_x \V).\D_{v}f
      -\gamma \D_{v}.\left(\D_{v}+v\right)f=0, \\
   & E(t,x) \defegal \D_x \Vnl(t,x)
   = -\zeta * \frac{\kappa}{|S^{d-1}|}
    \frac{x}{|x|^d}*_x \rho(t,x)  \ \
     \ \ \
     \text{where } \rho(t,x) = \int     f(t,x,v)dv, \\
    &   f|_{t=0} = f_0,
    \end{aligned}
    \right.
\end{equation}
where $\zeta \in \ss$ (depending only on $x$). We shall write in
the following
\begin{equation} \label{defphi} \phi= -\zeta *
\frac{1}{|S^{d-1}|}
    \frac{x}{|x|^d}
    \end{equation}
    so that the field  reads $E = \kappa \phi *\rho$.
    In fact the following result
     and the ones in the next section work as well for
     any $\phi \in L^\infty$.
       We shall prove the following Proposition.

\begin{prop} \label{existA}
Suppose that $V$ satisfies hypothesis \Hun and that  $f_0 \in
B^2$. Then for all $T>0$, the approximate problem
{\rm(\ref{VPFPA})} admits a unique strong  solution in $\cc([0,T[,
B^2) $.
\end{prop}

\preuve We suppose $\gamma=1$ ($\D_v= -b^*$) in the proof for
convenience. The existence is given by a fixed point theorem. We
study the following family of linear problems where $f^0$ is fixed
and on an interval of time $[0,T]$ for $T$ finite and fixed.
\begin{equation*}
  \label{VPFPAn}
  \left\{
  \begin{aligned}
   &  \D_{t}f^{n+1}+v.\D_{x}f^{n+1} -(E^n
   + \D_x \V).\D_{v}f^{n+1}
      -\gamma \D_{v}.\left(\D_{v}+v\right)f^{n+1}=0, \\
   & E^n= \D_x \Vnl^n = \kappa \phi*_x \rho^n   \ \
     \ \ \
     \text{with } \rho^n(t,x) = \int     f^n(t,x,v)dv, \\
    &     f|_{t=0} = f_0  \ \ \ \textrm{ and } \ \ \ f^0 = f_0.
    \end{aligned}
    \right.
\end{equation*}
In the following we call $C$ any constant independent of $n$ (but
perhaps depending on $T$). Proposition \ref{existmild} yields that
for each $n \geq 0$  this problem admits a mild solution $f^{n+1}
$ since
$$
\norm{E^n }_{L^{\infty}(dtdx)} \leq \norm{\phi}_{L^{\infty}}
\norm{\rho^n}_{L^1(dx)} =\norm{\phi}_{L^{\infty}}=C,
$$
from  Young inequality. This solution is given by
$$
f^{n+1}(t) = e^{-tK} f_0 - \int_0^t e^{-(t-s)K}b^* E^n
f^{n+1}(s)ds, ds
$$
and we observe using a Gronwall inequality and the diffusion
estimate from Proposition \ref{disp} that there exists a constant
$C$ independent of $n$ such that $\norm{f^n}_{B^2} \leq C_T$. Now
for all $0 \leq t \leq T$
\begin{equation*}
\begin{split}
\norm{f^{n+1} - f^n}_{B^2} & \leq \norm{ \int_0^t e^{-(t-s)K}b^*
E^n
(f^{n+1}(s) - f^n(s)) ds}_{B^2} \\
& \ \ \ \ + \norm{ \int_0^t e^{-(t-s)K}b^* f^n
(E^n(s) - E^{n-1}(s)) ds}_{B^2} \\
&  \leq C_T' \sqrt{t} (\norm{f^{n+1} - f^n}_{B^2} -
\norm{f^n-f^{n-1}}_{B^2})
\end{split}
\end{equation*}
since
$$
\norm{E^n(s) - E^{n-1}(s)}_{L^\infty} \leq
\norm{\phi}_{L^{\infty}} \norm{f^n-f^{n-1}}_{L^1} \leq
\norm{\phi}_{L^{\infty}}\norm{f^n-f^{n-1}}_{B^2}.
$$
Now a standard
fixed point theorem give that on any interval $[0, c_T[ \in [0,T[$
the scheme converges in $L^\infty( [0, c_T[, B^2)$ where $c_T$ is
independent of $n$. We can apply the same procedure on any
interval of type  $[t, t+c_T[ \subset [0,T[$ for $t$ arbitrary and
we get that $f^n$ converges (strongly) in $L^\infty([0,T[, B^2)$
toward a function $f$, and that this is also the case for $E^n$
toward $E$ in $L^\infty([0,T[, L^\infty)$ where $E$ is given by
$$
E= \kappa \phi *_x \int f(t,x,v) dv.
$$
The function $f$  is therefore a mild solution
 of the problem $\D_t f + K f + E b^* f = 0$, $f|_{t=0} = f_0$.
 Since by Proposition \ref{existmild} the solution is unique we
 get the result. \fin

\Subsection{Exponential time decay for small nonlinear coupling}

In this subsection we continue the study of the mollified
Vlasov-Poisson-Fokker-Planck equation defined in (\ref{VPFPA}).
 The
aim of this section is to prove Theorem \ref{maindecay} and
Corollary \ref{cordecay} about the exponential decay for small
charge.  Let us now define as in the Introduction the
Fokker-Planck operator corresponding to the stationnary
Vlasov-Poisson-Fokker Planck equation
\begin{equation*} \label{FPinf2}
K_\infty = v\D_x - \D_x(\V+\Vinf) \D_v -\gamma \D_v(\D_v +  v).
\end{equation*}
We know that $\Vinf \in W^{\infty,\infty}$ so that total potential
$\V+\Vinf$ satisfies Hypothesis (H1), and we suppose that it also
satisfies Hypothesis (H2bis). The Maxwellian associated to this
operator is
\begin{equation*}
\mminf(x,v) = \frac{e^{-(v^2/2 + \V(x)+\Vinf(x))}}{\int e^{-(v^2/2
+ \V(x)+\Vinf(x))}dxdv}
\end{equation*}
and is in $\ss \subset L^1$ with norm $1$ in $L^1$. We define also
the associated spaces $\Binf^2 = \set{ f \in \dd' \text{ s. t. }
f/\mminf \in L^2(\mminf dxdv)}$. Since $\V + \Vinf$ satisfy the
hypothesis (H1) and (H2bis) we can apply all the results obtained
for a generic potential $V$. We recall that $\ai$ is the smallest
positive real part  of the eigenvalues of the corresponding Witten
Laplacian
\begin{equation*}
\Lambda_\infty^2 =-\gamma \D_v(\D_v +  v) -\gamma \D_x(\D_x +
\D_x(\V +\Vinf))
\end{equation*}
in $\Binf^2$. We denote by the same symbols $\Lambda_\infty^2$ and
$\Kinf$ the closure from $\cc^\infty_0$ of the corresponding
operators in $\Binf^2$, and recall that they are maximal accretive
from Proposition \ref{cauchy}. We then follow Subsection
\ref{lindecay} by defining in our context the following space
$$
\Binfp^2 = \mm^{\perp} \cap \Binf^2 = \set{f \in \Binf^2 \text{
s.t. } \int f dxdv = 0},
$$
 endowed with the norm of $\Binf^2$,
  where $\perp$ stands for the orthogonal with respect to
the scalar product. We note that $\Binfp^2$ is stable for $\Kinf
$. The following proposition is a direct consequence of
Propositions \ref{disp2} and \ref{expdecay} for $\V + \Vinf$:

\begin{prop} \label{disp2close}
There exists  constants $C_\infty$ and $A_\infty$
 such that for all $t>0$, \\
 i)  $e^{-t\Kinf} b^* \textrm{ is bounded by } C_\infty
 (1+t^{-1/2})
e^{-\ai t/A_\infty}$ on $\Binf^2$ \\
ii) $ e^{-t\Kinf} a^* \textrm{ is bounded by } C_\infty
(1+t^{-3/2})
e^{-\ai t/A_\infty}$ on $\Binf^2$ \\
iii) $ e^{-t\Kinf}  \textrm{ is bounded by } 3  e^{-\ai
t/A_\infty}$ on
$\Binfp^2$ \\
 where $C_\infty$ and $A_\infty$ depend only
  $\norm{(\Ve + \Vinf)''}_{L^\infty}$ and the physical
   constants (uniformly in $\kappa$ varying in a fixed compact set).
\end{prop}

We work now in the Hilbert space $\Binf^2$ which we recall is
norm-equivalent to $B^2$ since $\Vinf \in W^{\infty,\infty}$. For
convenience we again suppose $\gamma =1$. For $t,x \in \R^+\times
\R^d$ we denote
$$
\Vd(t,x) = \Vnl(t,x) -\Vinf(x).
 $$
We can write the Cauchy problem associated to the VPFP system as
follows
\begin{equation*}
\left\{
\begin{aligned}
&  \D_t f + \Kinf f = -b^* \D_x \Vd f, \\
& f|_{t=0} = f_0.
\end{aligned}
\right.
\end{equation*}
Using the a priori bounds for the solution  $f$ given by
Proposition \ref{existA}, the unique solution satisfies  the
following Duhamel formula written in terms of $\Kinf$ in
$\Binf^2$:
\begin{equation} \label{duha1}
f(t,x,v) = e^{-t\Kinf} f_0 (x,v) - \int_0^t e^{-(t-s)\Kinf}
b^*\D_x \Vd (s,x) f(s,x,v)ds.
\end{equation}
We know that $ \D_x \Vd \in L^\infty(\R_t \times \R_x^d) $ and we
recall that $
 \phi(x) = -\zeta*\frac{1}{|S^{d-1}|}
\frac{x}{|x|^d},
$
 so that $\D_x \Vd$ reads
\begin{equation*}
\begin{split}
\D_x \Vd(t,x) & = \D_x \Vnl(t,x) - \D_x V_{\infty}(x) = \kappa
\phi(x) *_x (\rho(t,x) - \rho_\infty(x))
\end{split}
\end{equation*}
where $\rho_\infty(t,x) = \int f_\infty(t,x,v) dv$ with
$f_\infty=\mminf$.
 This is clear that $f_\infty=\mminf$ is the projection in the
 Hilbert space $\Binf^2$
of the Cauchy data $f_0$ on the fundamental space
$\text{Span}(\mminf)$ since
\begin{equation*} \label{finf}
f_\infty = (f_0,\mminf)_{\Binf^2} \mminf = \sep{\iint f_0 \mminf
\mminf^{-1}dxdv} \mminf = \mminf.
\end{equation*}
 Let us denote
by $g(t,x,v) = f(t,x,v) - f_\infty(x,v)$. Since $f_\infty \in
\text{Ker}( K_\infty)$ we have
$$
e^{-t\Kinf}f_\infty = f_\infty.
$$
The Duhamel formula (\ref{duha1}) therefore reads
$$
g(t)  =
e^{-t\Kinf} g_0 + \int_0^t e^{-(t-s)\Kinf} b^*\D_x \Vd (s) f(s)ds
$$
and we have $ \D_x \Vd(t,x) = \kappa \phi *_x \int g(t,x,v) dv. $
In fact we shall use the following representation
\begin{equation} \label{duhamel}
 g(t)  =   e^{-t\Kinf} g_0 + \int_0^t e^{-(t-s)\Kinf} b^*\D_x \Vd
(s) (g(s)+f_\infty) ds
\end{equation}

Now we  take the $\Binf^2$ norm in this formula. We first note
that $g_0 \in \Binfp^2$ which gives from Proposition
\ref{disp2close} that  for all $t \geq 0$,
\begin{equation} \label{bp1}
\norm{e^{-t\Kinf} g_0}_{\Binf^2} \leq 3 e^{-\ai t /\Ai}
\norm{g_0}_{\Binf^2}
\end{equation}
In order to estimate the integrals in (\ref{duhamel}), we first
estimate the $L^\infty(dx) $ norm of $\D_x \Vd(s)$ for all $s \in
(0,t)$. In the following $\norm{\phi}$ stands for
$\norm{\phi}_{L^\infty}$. First we note that
$$
\norm{ \D_x \Vd(s)}_{L^{\infty}(dx)}  \leq   \kappa \norm{\phi}
\norm{ \int g(s) dv }_{L^1(dx)} = \kappa \norm{\phi} \norm{g(s)
}_{L^1(dxdv)} \leq \kappa \norm{\phi}\norm{g(s) }_{\Binf^2}
$$
which gives
\begin{equation} \label{suite1}
\begin{split}
\norm{\D_x \Vd (s) f_\infty(s) }_{\Binf^2} & \leq  \norm{\D_x
\Vd(t)}_{L^\infty} \norm{f_\infty}_{\Binf^2}
 \leq   \kappa \norm{\phi} \norm{g(s)}_{\Binf^2}
\end{split}
\end{equation}
Now we estimate the norm of $\D_x \Vd(s)$ in an another way
\begin{equation*}
\begin{split}
\norm{ \D_x \Vd(s)}_{L^{\infty}} & \leq  \kappa \norm{\phi}
\norm{g(s) }_{L^1(dxdv)} \leq \kappa \norm{\phi} (\norm{f(s)
}_{L^1(dxdv)} + \norm{f_\infty }_{L^1(dxdv)} )\leq  2\kappa
\norm{\phi}
\end{split}
\end{equation*}
since $f$ and $f_\infty$ are $L^1$ normalized. This gives
\begin{equation} \label{suite2}
\begin{split}
\norm{\D_x \Vd (s) g(s) }_{\Binf^2} & \leq  \norm{\D_x
\Vd(s)}_{L^\infty} \norm{g(s)}_{\Binf^2}
 \leq   2\kappa \norm{\phi} \norm{g}_{\Binf^2}.
\end{split}
\end{equation}
Putting together (\ref{suite1}-\ref{suite2}) we get
\begin{equation*} \label{suite3}
\begin{split}
\norm{\D_x \Vd (s) (g(s)+f_\infty) }_{\Binf^2} &\leq 3 \kappa
\norm{\phi} \norm{g}_{\Binf^2}.
\end{split}
\end{equation*}
 Now applying  Proposition
\ref{disp2close} to the operator $\Kinf$ with the associated rate
$\ai$. We can write for $t-s > 0$
\begin{equation} \label{bp2}
\begin{split}
& \norm{e^{-(t-s)\Kinf} b^*\D_x \Vd (s) (g(s)+f_\infty)}_{\Binf^2}  \\
& \leq C_2 (1+(t-s)^{-1/2}) e^{-\ai (t-s)/\Ai}
\norm{\D_x \Vd (s) (g(s)+f_\infty)}_{\Binf^2} \\
& \leq  \kappa C(1+(t-s)^{-1/2}) e^{-\ai
(t-s)/\Ai}\norm{g}_{\Binf^2}
\end{split}
\end{equation}
Putting (\ref{bp1}-\ref{bp2}) in the Duhamel Formula
(\ref{duhamel}) and calling from now on $C$ any  constant
depending on $\norm{\phi}$ and the derivatives of $V+\Vinf$, we
get
\begin{equation*} \label{duhamel2}
\begin{split}
 \norm{g(t)}_{\Binf^2} \leq & 3 e^{-\ai t/\Ai}
\norm{g_0}_{\Binf^2} \\
& + C \kappa \int_0^t (1+(t-s)^{-1/2}) e^{-\ai (t-s)/\Ai}
 \norm{g(s)}_{\Binf^2} ds
\end{split}
\end{equation*}
Let us define for $t\geq 0$,  $\psi(t) = e^{\ai t/(2\Ai) }
 \norm{g(t)}_{\Binf^2} $.  We get for $t\geq 0$,
$$
\psi(t) \leq 3 \psi(0) + C \kappa  \int_0^t (1+(t-s)^{-1/2})
e^{-\ai (t-s)/(2\Ai)}  \psi(s)  ds .
$$
With an other constant $C$ we get
$$
\psi(t) \leq 3 \psi(0) + ( C \kappa / \ai )
 \sup_{s \in [0,t]}  \psi(s).
$$
Note here that contrary to $\Ai$ the constant $\ai$ cannot be
absorbed in the constant $C$ since not controlled by semi-norms of
$(\Ve + \Vinf)''$. Under the following assumption
$$
C \kappa /\ai \leq 1/2
$$
we get that for all $t\geq 0$, $  \sup_{s \in [0,t]}  \psi(s) \leq
6 \psi(0). $ This reads in terms of $g$:
$$
 \norm{g(t)}_{\Binf^2} \leq 6
 \norm{g(0)}_{\Binf^2} e^{-\ai t/(2\Ai) },
 $$
 and
 the proof of Theorem \ref{maindecay} is complete. \fin

\preuve[of Corollary \ref{cordecay}] First recall that $ 0 \leq
H({f}(t)|\mminf) $ since $f$ and $\mminf$ are $L^1$-normalized.
Using the inequality $\ln(s) \leq s-1$ we get
\begin{equation*}
\begin{split} H({f}(t)|\mminf) & = \iint
\frac{{f}(t)}{\mminf} \ln \sep{\frac{{f}(t)}{\mminf}}
\mminf dxdv \\
& \leq \iint \frac{{f}(t)}{\mminf} \sep{\frac{{f}(t)
-\mminf}{\mminf}} \mminf dxdv \\
& \leq  \norm{{f}(t)}_{\Binf^2} \norm{{f}(t)- \mminf}_{\Binf^2}
\end{split}
\end{equation*}
Now applying Theorem \ref{maindecay} we first notice that
$$
\norm{{f}(t)}_{\Binf^2} \leq
\norm{{g}(t)}_{\Binf^2}+\norm{f_\infty}_{\Binf^2} \leq 7
\norm{f_0}_{\Binf^2}
$$
since $\norm{f_\infty}_{\Binf^2}=1 \leq \norm{f_0}_{\Binf^2}$.
Using again Theorem \ref{maindecay} yields
$$
H(f(t)|f_\infty) \leq C
\norm{f_0}_{\Binf^2}\norm{f_0-f_\infty}_{\Binf^2} e^{-\ai t
/2\Ai}.
$$
The proof of Corollary \ref{cordecay} is complete. \fin


\end{document}